# New conjectures in number theory - The distribution of prime numbers -

Jonas Castillo Toloza

*"Mathematicians have tried in vain to this day to discover some order in the sequence of prime numbers, and we have reason to believe that it is a mystery into which the human mind will never penetrate". (LEONHARD EULER).*

**Abstract**

Since the mathematicians of ancient Greece until Fermat, since Gauss until today; the way how the primes along the numerical straight line are distributed has become perhaps the most difficult math problem; many people believe that their distribution is chaotic, governed only by the laws of chance. In this article the author tries to resolve satisfactorily all matter related to the distance between primes, based on the determination of a lower bound for the primes counting function, the estimate of Legendre's conjecture, and the iterative method for solving equations.

*To Elian David, the locomotive of my dreams and to the memory of the eminent indian mathematician Srinivasa Ramanujan Aiyangar, master of all self-taught.*

## 1. The lower bound for the primes counting function.

If $\pi(N)$ is the number of primes less than or equal to $N$, and $N \geq 4$, then

$$\frac{N}{\log_2 N}$$

is the lower bound for $\pi(N)$.

Understood here that $\log_2 N$ is the logarithm base $2$ of $N$.

The prime number theorem states that

$$\frac{N}{\ln N} \sim \pi(N)$$

Rosser and Schoenfeld proved that if $N \geq 17$, then

$$\pi(N) > \frac{N}{\ln N}$$

Moreover, a result of Erdős-Kalmár affirms that

Given

$$\varepsilon > 0 \text{ , } \exists N_0 = N_0(\varepsilon)$$

such that

$$(\log 2 - \varepsilon)\frac{N}{\log N} \leq \pi(N) \leq (\log 4 - \varepsilon)\frac{N}{\log N} \text{ , for } N \geq N_0$$

Where it is deduced by the first inequality that

$$\frac{N}{\log_2 N}$$

Is the lower bound for $\pi(N)$.

Henceforth whenever we find the lower bound of primes located in a specific range, where $\ln N$ apperars in the formulas, we replace it with $\log_2 N$ .

## 2. Conjecture I. The estimate of Legendre's conjecture.

If $\pi(N,M)$ is the number of prime numbers ranging from $N$ to $M$ excluding extremes, then

$$\pi\left(a^2,(a+1)^2\right) \approx \frac{1}{2}\pi(2a+1)$$

The Legendre's conjecture states that there always exist at least one prime between $a^2$ and $(a+1)^2$ . If we accept that Legendre's Conjecture is true, then there should be a pattern to measure the number of prime numbers between a square number and the next, finding the author that “the number of primes located between $a^2$ and $(a+1)^2$ is approximately equal to half of prime numbers that are less than or equal to $(2a+1)$“. Through the prime number theorem we express as follows:

$$\pi\left(a^2,(a+1)^2\right) \approx \frac{2a+1}{2\ln(2a+1)}$$

The empirical evidence in favor of this conjecture is enormous.

*Table 1. The estimate of Legendre's conjecture.*

| $a$ | $\pi\left(a^2,(a+1)^2\right)$ | $\frac{1}{2}\pi(2a+1)$ |
|---|---|---|
| 10 | 5 | 4 |
| 20 | 7 | 6 |
| $\sqrt{1000}$ | 11 | 9 |
| 100 | 23 | 23 |
| 213,48 | 37 | 41 |
| 500 | 71 | 84 |
| 1 000 | 152 | 151 |
| 2 000 | 267 | 275 |
| 5 000 | 613 | 614 |
| 10 000 | 1 081 | 1 131 |
| 20 000 | 2 020 | 2 101 |
| 50 000 | 4 605 | 4 796 |
| 100 000 | 8 668 | 8 992 |
| 200 000 | 16 473 | 16 930 |
| 300 000 | 23 965 | 24 549 |
| 500 000 | 38 250 | 39 249 |
| 800 000 | 59 091 | 60 563 |
| 1 000 000 | 72 413 | 74 466 |

## 2.1. Proposition.

If $a \geq 1/2$, the lower bound for $\pi\left(a^2,(a+1)^2\right)$ is

$$\left\lfloor \frac{2a+1}{2\log_2(2a+1)} \right\rfloor$$

Let's take the example when $a = 14.2$, we have the lower bound for $\pi\left(14.2^2, 15.2^2\right)$ is

$$\left\lfloor \frac{29.4}{2\log_2 29.4} \right\rfloor = 3$$

## 3. Conjecture II.

If $\sqrt{b} - \sqrt{a} \geq 1$ (taking as variables $a$ and $b$), then

$$\pi(a,b) \approx \frac{b-a}{2\ln\left(\sqrt{b}+\sqrt{a}\right)}$$

Example. $a = 70$, $b = 100$

We verify

$$\sqrt{100} - \sqrt{70} > 1$$

holds, then

$$\pi(70,100) \approx \frac{30}{2\ln\left(\sqrt{100}+\sqrt{70}\right)} = 5{,}15369...$$

There are 6 prime numbers in that range.

Let's look at case when $a = 0$, we have that

$$\pi(0,b) = \pi(b) \approx \frac{b-0}{2\ln\left(\sqrt{b}+\sqrt{0}\right)} \approx \frac{b}{\ln b}$$

According to the stament of the prime number theorem.

Consider another special case

When

$$s \to 2(+), \; a = \left(\sqrt{s}+1\right)^{\frac{2}{s-2}}$$

then

$$\pi(a, a.s) \approx \pi(a)$$

"The number of prime numbers ranging from $a$ to $a.s$ approximates to the number of prime numbers less than or equal to $a$".

Obtaining the result.

Doing $b = a.s$ we have the conjecture II that

$$(1)\quad \pi(a, a.s) \approx \frac{a.s - a}{2\ln\left(\sqrt{a.s} + \sqrt{a}\right)}$$

By the prime number theorem we have

$$(2)\quad \pi(a) \approx \frac{a}{\ln a}$$

equating (1) and (2) we have

$$\frac{a.s - a}{2\ln\left(\sqrt{a.s} + \sqrt{a}\right)} \approx \frac{a}{\ln a}$$

$$\frac{a(s-1)}{2\ln\left(\sqrt{a}\left(\sqrt{s}+1\right)\right)} \approx \frac{a}{\ln a}$$

$$\frac{s-1}{2\ln\left(\sqrt{a}\left(\sqrt{s}+1\right)\right)} \approx \frac{1}{\ln a}$$

$$\frac{s-1}{2\ln\sqrt{a} + 2\ln\left(\sqrt{s}+1\right)} \approx \frac{1}{\ln a}$$

$$\frac{s-1}{\ln a + 2\ln\left(\sqrt{s}+1\right)} \approx \frac{1}{\ln a}$$

$$\ln a.(s-1) \approx \ln a + 2\ln\left(\sqrt{s}+1\right)$$

$$\ln a.(s-1) - \ln a \approx 2\ln\left(\sqrt{s}+1\right)$$

$$\ln a.(s-2) \approx 2\ln\left(\sqrt{s}+1\right)$$

$$\ln a.(s-2) \approx \ln\left(\sqrt{s}+1\right)^2$$

$$\ln a \approx \frac{\ln\left(\sqrt{s}+1\right)^2}{s-2}$$

$$a \approx \left(\sqrt{s}+1\right)^{\frac{2}{s-2}}$$

Solution of the equation

$$a = \left(\sqrt{s}+1\right)^{\frac{2}{s-2}}$$

$$s = \frac{t}{a}+1$$

We find the value of $t$ by the iterative method

$$t_0 = \frac{2a.\ln\left(2\sqrt{a}+1\right)}{\ln a}$$

$$t_1 = \frac{2a.\ln\left(2\sqrt{a+t_0}+\sqrt{a}\right)}{\ln a}$$

$$t_2 = \frac{2a.\ln\left(2\sqrt{a+t_1}+\sqrt{a}\right)}{\ln a}$$

$$t_3 = \frac{2a.\ln\left(2\sqrt{a+t_2}+\sqrt{a}\right)}{\ln a}$$

$$t_m = \frac{2a.\ln\left(2\sqrt{a+t_{(m-1)}}+\sqrt{a}\right)}{\ln a}$$

When $m \to \infty$ then $t_m \to t$

A well known result in analytic number theory says that if

$x(s)$ and $y(s)$ are funtions of the parameter $s$ such that $x(s)$ tends to the infinity with $s$ and limit of $y(s)/x(s)=2$, then $\pi\left(x(s),y(s)\right)$ it is asymptotically $\pi\left(x(s)\right)$ .

Here we come to the same result, but for a totally different way, being in harmony conjecture II with known results.

## 3.1. Proposition

If $4 \le b-a$, $\sqrt{b}-\sqrt{a} \ge 1$, the lower bound for $\pi(a,b)$ is

$$\left\lfloor \frac{b-a}{2\log_2\left(\sqrt{b}+\sqrt{a}\right)} \right\rfloor$$

Example $a=200$, $b=230$

We verify

$$4<30,\ \sqrt{230}-\sqrt{200}>1$$

holds, then the lower bound for $\pi(200,230)$ is

$$\left\lfloor \frac{30}{2\log_2\left(\sqrt{230}+\sqrt{200}\right)} \right\rfloor = 3$$

This result encompasses Brocard's Conjecture, which says that there are at least four primes between $(p_n)^2$ and $(p_{n+1})^2$, being $p_n$ the nth prime number. We have the shortest distance between two consecutive prime is 2, and for this proposition we deduce that if $a > 10.23$ then the lower bound for $\pi\left(a^2,(a+2)^2\right)$ is greater than $4$, regardless of $a$ is a prime number or it is not, even $a$ can be non-integer number.

Example. If $a = 127$ we have between $127^2$ and $129^2$ there are at least $32$ primes. But the proposition's conjecture I it is very strong and we have that if $a \geq 28.89$ then the lower bound for

$\pi\left(a^2,(a+1)^2\right)$ is greater than $4$.

## 4. conjecture III.

If $2 \leq b-a$, $\sqrt{b}-\sqrt{a} \leq 1$, the lower bound for $\pi(a,b)$ is

$$\left\lfloor \frac{(b-a)\left(\sqrt{b}-\sqrt{a}\right)^2}{2\log_2(b-a)} \right\rfloor$$

Example. $a = 113.3$, $b = 129$

We verify

$$2 < 15.7,\ \sqrt{129}-\sqrt{113.3} < 1$$

holds, then the lower bound for $\pi(113.3,129)$ is

$$\left\lfloor \frac{15.7\left(\sqrt{129}-\sqrt{113.3}\right)^2}{2\log_2 15.7} \right\rfloor = 1$$

## 5. Development of algorithms.

If $n$ is the lower bound for $\pi(a,b)$ and $x = (b-a)$, and if the conjectures II and III are true, we can develop algorithms to calculate the values of $a$, $b$ and $x$; being $x$ the maximum distance between $a$ and the nth prime number that happens, or also the maximum distance between $b$ and the nth prime number preceding it; we mean that however great the distance between $a$ and the nth prime number that follow, this is never over to $x$.

The iterative method for solving equations plays an important role in the development of these algorithms.

## 5.1. Algorithms for the calculation of $x$ given the values of $a$ and $n$.

**5.1.1.** If

$$\frac{2\sqrt{a}+1}{2\log_2\left(2\sqrt{a}+1\right)}=n$$

then

$$x=2\sqrt{a}+1$$

We take $x=\lceil 2\sqrt{a}+1\rceil$

Example. Find the value of $x$ when $a=16256.25$, $n=16$.

We verify

$$\frac{2\sqrt{16256.25}+1}{2\log_2\left(2\sqrt{16256.25}+1\right)}=16$$

holds, then

$$x=2\sqrt{16256.25}+1=256$$

We have that the maximun distance betwen $16256.25$ and the 16th prime number that follows it is $256$. However great the distance betwen $16256.25$ and the 16th prime number that follows, it never exceeds $256$.

By the proposition of the conjecture II or through the conjecture III, we have that the lower bound for $\pi\left(16256.25,16512.25\right)$ is $16$.

**5.1.2.** If

$$\frac{2\sqrt{a}+1}{2\log_2\left(2\sqrt{a}+1\right)}\geq n$$

Then we calculate the value of $x$ by the iterative method.

$$x_0=\sqrt{2n\left(2\sqrt{a}+1\right).\log_2\left(2\sqrt{a}+1\right)}$$

$$x_1 = \frac{\left(\sqrt{a + x_0} + \sqrt{a}\right)\sqrt{2n.x_0.\log_2 x_0}}{x_0}$$

$$x_2 = \frac{\left(\sqrt{a + x_1} + \sqrt{a}\right)\sqrt{2n.x_1.\log_2 x_1}}{x_1}$$

$$x_3 = \frac{\left(\sqrt{a + x_2} + \sqrt{a}\right)\sqrt{2n.x_2.\log_2 x_2}}{x_2}$$

$$x_m = \frac{\left(\sqrt{a + x_{m-1}} + \sqrt{a}\right)\sqrt{2n.x_{m-1}.\log_2 x_{m-1}}}{x_{m-1}}$$

when $m \to \infty$ then $x_m \to x$

We take $x = \lceil x_m \rceil$

Example. Find the value of $x$ when $a = 1024$, $n = 2$.

We verify

$$\frac{2\sqrt{1024} + 1}{2\log_2\left(2\sqrt{1024} + 1\right)} > 2$$

holds, then we calculate the value of $x$ by the iterative method.

$$x_0 = \sqrt{4 * 65 * \log_2 65} = 39.570388...$$

$$x_1 = \frac{\left(\sqrt{1024 + x_0} + \sqrt{1024}\right)\sqrt{4x_0 * \log_2 x_0}}{x_0} = 47.321533...$$

$$x_2 = 44.3939...$$

$$x_3 = 45.4218...$$

$$x_4 = 45.05119...$$

We take $x = 46$.

So we have that the maximun distance betwen 1024 and the second prime number that follow it is $46$. We verify through the conjecture III. The lower bound for $\pi(1024, 1070)$ is

$$\left\lfloor \frac{46\left(\sqrt{1070} - \sqrt{1024}\right)^2}{2\log_2 46} \right\rfloor = 2$$

**5.1.3.** If

$$\frac{2\sqrt{a}+1}{2\log_2\left(2\sqrt{a}+1\right)} \le n$$

Then we calculate the value of $x$ by the iterative method.

$$x_0 = 2n.\log_2\left(2\sqrt{a}+1\right)$$

$$x_1 = 2n.\log_2\left(2\sqrt{a+x_0}+\sqrt{a}\right)$$

$$x_2 = 2n.\log_2\left(2\sqrt{a+x_1}+\sqrt{a}\right)$$

$$x_3 = 2n.\log_2\left(2\sqrt{a+x_2}+\sqrt{a}\right)$$

$$x_m = 2n.\log_2\left(2\sqrt{a+x_{m-1}}+\sqrt{a}\right)$$

when $m \to \infty$ then $x_m \to x$

we take $x = \lceil x_m \rceil$

Example. Find the value of $x$ when $a = 10000,\ n = 36$.

We verify

$$\frac{2\sqrt{10000}+1}{2\log_2\left(2\sqrt{10000}+1\right)} < 36$$

holds, then we calculate the value of $x$ by the iterative method

$$x_0 = 72\log_2 201 = 550.8757...$$

$$x_1 = 72\log_2\left(2\sqrt{10000+x_0}+\sqrt{10000}\right) = 551.759...$$

$$x_2 = 551.761...$$

We take $x = 552$.

So we have the maximu distance between $10000$ and the 36th prime number that follows it is $552$. We verify by the proposition of the conjecture II. The lower bound for $\pi(10000,10552)$ is

$$\left\lfloor \frac{552}{2\log_2\left(\sqrt{10552}+\sqrt{10000}\right)} \right\rfloor = 36$$

## 5.2. Algorithms for the calculation of $x$ given the values of $b$ and $n$.

**5.2.1.** If

$$\frac{2\sqrt{b}-1}{2\log_2\left(2\sqrt{b}-1\right)}=n$$

Then

$$x=2\sqrt{b}-1$$

We take $x=\lceil 2\sqrt{b}-1\rceil$

Example. Find the value of $x$ when $b=16512.25$, $n=16$.

We verify

$$\frac{2\sqrt{16512.25}-1}{2\log_2\left(2\sqrt{16512.25}-1\right)}=16$$

Then

$$x=2\sqrt{16512.25}-1=256$$

We have the maximun distance betwen $16512.25$ and the 16th prime number that precedes it is $256$, according to the solution of the example of **5.1.1**.

**5.2.2.** If

$$\frac{2\sqrt{b}-1}{2\log_2\left(2\sqrt{b}-1\right)}\geq n$$

Then we calculate the value of $x$ by the iterative method.

$$x_0=\sqrt{2n\left(2\sqrt{b}-1\right).\log_2\left(2\sqrt{b}-1\right)}$$

$$x_1=\frac{\left(\sqrt{b-x_0}+\sqrt{b}\right)\sqrt{2n.x_0.\log_2 x_0}}{x_0}$$

$$x_2=\frac{\left(\sqrt{b-x_1}+\sqrt{b}\right)\sqrt{2n.x_1.\log_2 x_1}}{x_1}$$

$$x_3=\frac{\left(\sqrt{b-x_2}+\sqrt{b}\right)\sqrt{2n.x_2.\log_2 x_2}}{x_2}$$

$$x_m = \frac{\left(\sqrt{b - x_{m-1}} + \sqrt{b}\right)\sqrt{2n.x_{m-1}.\log_2 x_{m-1}}}{x_{m-1}}$$

when $m \to \infty$ then $x_m \to x$

We take it $x = \lceil x_m \rceil$

Example. Find the value of $x$ when $b = 1070$, $n = 2$.

We verify

$$\frac{2\sqrt{1070} - 1}{2\log_2\left(2\sqrt{1070} - 1\right)} > 2$$

holds, then we calculate the value of $x$ by the iterative method

$$x_0 = \sqrt{4\left(2\sqrt{1070} - 1\right) * \log_2\left(2\sqrt{1070} - 1\right)} = 39.3517...$$

$$x_1 = \frac{\left(\sqrt{1070 - x_0} + \sqrt{1070}\right)\sqrt{4x_0.\log_2 x_0}}{x_0} = 45.565...$$

$$x_2 = 44.2788...$$

$$x_3 = 45.5014...$$

$$x_4 = 45.0536...$$

We take $x = 46$

We have that the maximun distance between $1070$ and the second prime number that precedes it is $46$, according with the solution of the example of **5.1.2.**

**5.2.3.** If

$$\frac{2\sqrt{b} - 1}{2\log_2\left(2\sqrt{b} - 1\right)} \le n$$

Then we calculate the value of $x$ by the iterative method.

$$x_0 = 2n.\log_2\left(2\sqrt{b} - 1\right)$$

$$x_1 = 2n.\log_2\left(2\sqrt{b - x_0} + \sqrt{b}\right)$$

$$x_2 = 2n.\log_2\left(2\sqrt{b - x_1} + \sqrt{b}\right)$$

$$x_3 = 2n.\log_2\left(2\sqrt{b - x_2} + \sqrt{b}\right)$$

$$x_m = 2n.\log_2\left(2\sqrt{b-x_{m-1}}+\sqrt{b}\right)$$

when $m \to \infty$ so $x_m \to x$

we take $x = \lceil x_m \rceil$

Example. Find the value of $x$ when $b = 10552$, $n = 36$.

We verify

$$\frac{2\sqrt{10552}-1}{2\log_2\left(2\sqrt{10552}-1\right)} < 36$$

holds, then we calculate the value of $x$ by the iterative method.

$$x_0 = 72\log_2\left(2\sqrt{10552}-1\right) = 552.641...$$

$$x_1 = 72\log_2\left(2\sqrt{10552-x_0}+\sqrt{10552}\right) = 551.760...$$

$$x_2 = 551.762...$$

We take $x = 552$.

We have the maximun distance between $10552$ and the 36th prime number that precedes is prime number that above is $552$, according with tha solution of example **5.1.3.**

## 5.3. Algorithms for the calculation of $a$ given the values of $x$ and $n$.

**5.3.1.** If

$$\frac{x}{2\log_2 x} \geq n$$

Then

$$a = \frac{\left(x^2 - 2n.\log_2 x\right)^2}{8n.x.\log_2 x}$$

Obtaining the result.

By the conjecture III we have that if $2 \leq b-a$, $\sqrt{b}-\sqrt{a} \leq 1$, then the lower bound for $\pi(a,b)$ is

$$\left\lfloor \frac{(b-a)\left(\sqrt{b}-\sqrt{a}\right)^2}{2\log_2(b-a)} \right\rfloor$$

Doing $b=(a+x)$ we have that the lower bound for $\pi(a,(a+x))$ is:

$$\left\lfloor \frac{x\left(\sqrt{a+x}-\sqrt{a}\right)^2}{2\log_2 x} \right\rfloor$$

We take

$$n=\frac{x\left(\sqrt{a+x}-\sqrt{a}\right)^2}{2\log_2 x}$$

We develop the equation to solve for $a$

$$\frac{2n.\log_2 x}{x}=\left(\sqrt{a+x}-\sqrt{a}\right)^2$$

$$\frac{2n.\log_2 x}{x}=2a+x-2\sqrt{a\left(a+x\right)}$$

$$2\sqrt{a\left(a+x\right)}=2a+x-\frac{2n.\log_2 x}{x}$$

$$2\sqrt{a(a+x)}=2a+\frac{x^2-2n.\log_2 x}{x}$$

Squaring both sides of the equation, we have

$$4a.\left(a+x\right)=4a^2+\frac{4a.\left(x^2-2n.\log_2 x\right)}{x}+\frac{\left(x^2-2n.\log_2 x\right)^2}{x^2}$$

$$4a.\left(a+x\right)-4a^2-\frac{4a.\left(x^2-2n.\log_2 x\right)}{x}=\frac{\left(x^2-2n.\log_2 x\right)^2}{x^2}$$

$$4a.\left(a+x-a-\frac{\left(x^2-2n.\log_2 x\right)}{x}\right)=\frac{\left(x^2-2n.\log_2 x\right)^2}{x^2}$$

$$4a.\left(x-\frac{\left(x^2-2n.\log_2 x\right)}{x}\right)=\frac{\left(x^2-2n.\log_2 x\right)^2}{x^2}$$

$$4a.\frac{\left(x^2-x^2+2n.\log_2 x\right)}{x}=\frac{\left(x^2-2n.\log_2 x\right)^2}{x^2}$$

$$\frac{8a.n.\log_2 x}{x}=\frac{\left(x^2-2n.\log_2 x\right)^2}{x^2}$$

$$\frac{8a.n.x.\log_2 x}{x^2} = \frac{\left(x^2 - 2n.\log_2 x\right)^2}{x^2}$$

We eliminate denominators

$$8a.n.x.\log_2 x = \left(x^2 - 2n.\log_2 x\right)^2$$

Then

$$a = \frac{\left(x^2 - 2n.\log_2 x\right)^2}{8n.x.\log_2 x}$$

Example. Find the value of $a$ when $x = 200$, $n = 3$

We verify

$$\frac{200}{2\log_2 200} > 3$$

holds, then

$$a = \frac{\left(200^2 - 6\log_2 200\right)^2}{24 * 200 * \log_2 200} = 43508$$

We have that between $43508$ and $43708$ there are at least $3$ prime numbers.

We verify by the conjecture III .The lower bound for $\pi\left(43508, 43708\right)$ is

$$\left\lfloor \frac{200\left(\sqrt{43708} - \sqrt{43508}\right)^2}{2\log_2 200} \right\rfloor = 3$$

It is noteworthy that

$$\left(a + x\right) = \frac{\left(x^2 + 2n.\log_2 x\right)^2}{8n.x.\log_2 x}$$

We therefore conclude that between

$$\frac{\left(x^2 - 2n.\log_2 x\right)^2}{8n.x.\log_2 x} \quad \text{and} \quad \frac{\left(x^2 + 2n.\log_2 x\right)^2}{8n.x.\log_2 x}$$

There are at least $n$ prime numbers, for

$$\frac{x}{2\log_2 x} \geq n$$

THE CASE $n = 1$

Between

$$\frac{\left(x^2 - 2\log_2 x\right)^2}{8x.\log_2 x} \quad \text{and} \quad \frac{\left(x^2 + 2\log_2 x\right)^2}{8x.\log_2 x}$$

there are at least 1 prime number, for $x \geq 4$

Result is better understood if we plot the functions $f(x)$ and $g(x)$ on the same graph.

$$f(x) = \frac{\left(x^2 - 2\log_2 x\right)^2}{8x.\log_2 x}, \quad g(x) = \frac{\left(x^2 + 2\log_2 x\right)^2}{8x.\log_2 x}$$

We have that between $f(x)$ and $g(x)$ there always exist at least one prime number.

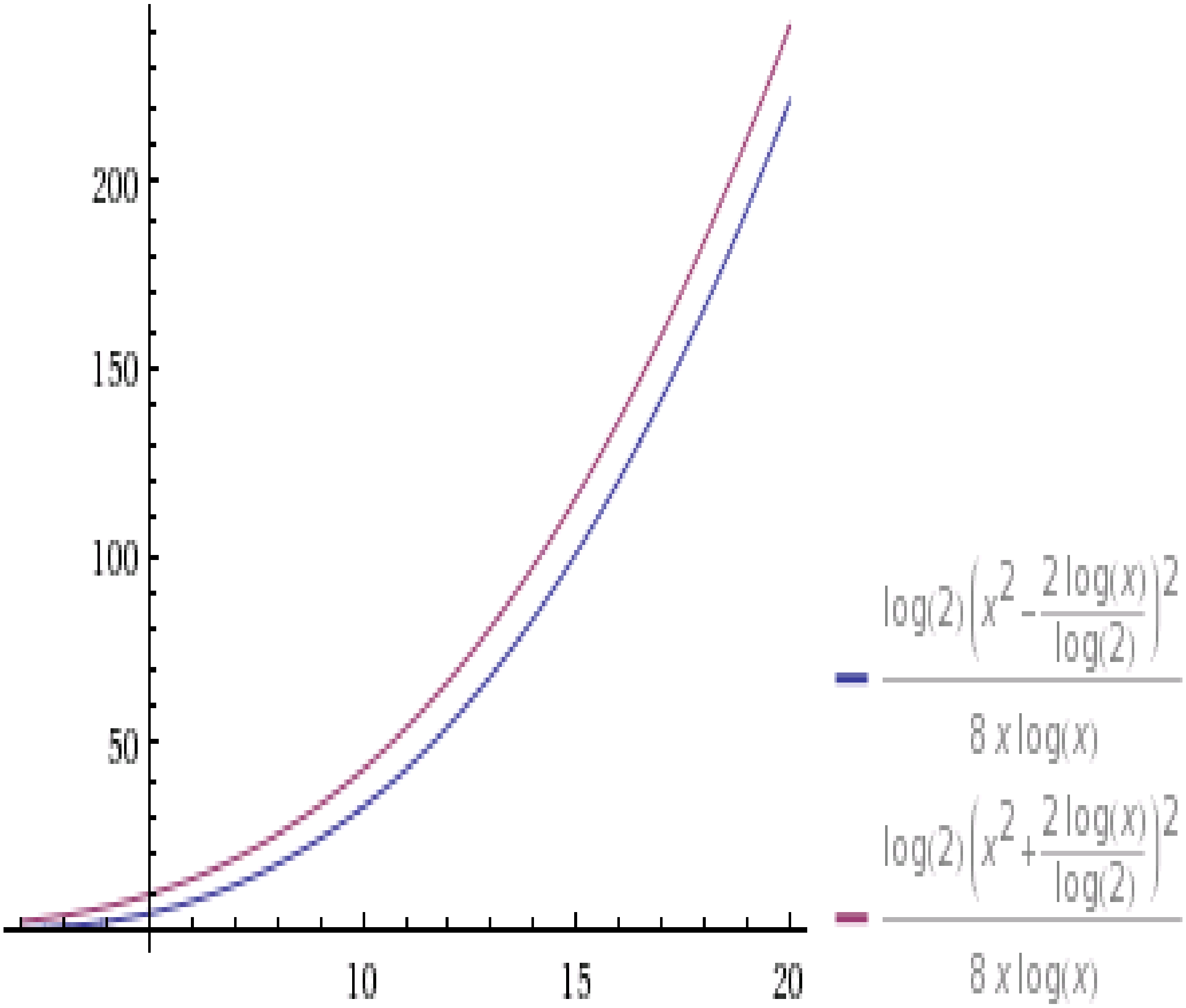

Graph taken from http://www.wolframalpha.com/input/?i=plot+%28%28%28x^2%29+-2log[2%2C%E2%81%A1x]+%29^2%29+%2F%288xlog[2%2C%E2%81%A1x]+%29%2C%28%28%28x^

2%29%2B2log%E2%81%A1[2%2Cx+]%29^2%29%2F%288xlog%E2%81%A1[2%2Cx+]%29+for+x%3D +2+to+20

The Andrica's conjecture says that

$$\sqrt{p_{n+1}}-\sqrt{p_n}<1$$

Being $p_n$ the nth prime number.

If we do $p_n=P$, $p_{n+1}=(P+x)$ so we have by **5.3.1**. that

$$\sqrt{P+x}-\sqrt{P}<\sqrt{\frac{2\log_2 x}{x}}$$

Taking here a criterion to determine when it is possible that $P$ and $(P+x)$ been consecutive prime numbers, for everything $P>2$

Note that

$$\sqrt{\frac{2\log_2 x}{x}}<1$$

For $x>4$ . We could say it:

$$\sqrt{p_{n+1}}-\sqrt{p_n}<\sqrt{\frac{2\log_2(p_{n+1}-p_n)}{p_{n+1}-p_n}}$$

For $n>1$. Conjecture Andrica reducing to a special exercise.

Take for example the numbers 199 and 223, and we want to determinate if it posible that those numbers are consecutive prime numbers.

We verify

$$\sqrt{223}-\sqrt{199}<\sqrt{\frac{2\log_2 14}{14}}$$

It is not met, then it is not possible that 199 and 223 are consecutive primes and therefore there is at least one prime number in that interval. The proof of the conjecture III (from which it follows **5.3.1**.) would demonstrate the Andrica's conjecture.

**5.3.2.** If

$$\frac{x}{\log_2 x}\geq n\geq\frac{x}{2\log_2 x}$$

Then

$$a = \frac{\left(2^{x/n} - x\right)^2}{4 * 2^{x/n}}$$

Obtaining the result.

By the proposition of the conjecture II we have that if $4 \le b - a$, $\sqrt{b} - \sqrt{a} \ge 1$, then the lower bound for $\pi(a, b)$ is

$$\left\lfloor \frac{b - a}{2\log_2\left(\sqrt{b} + \sqrt{a}\right)} \right\rfloor$$

Doinng $b = (a + x)$ we have that the lower bound for $\pi\left(a, (a + x)\right)$ is

$$\left\lfloor \frac{x}{2\log_2\left(\sqrt{a + x} + \sqrt{a}\right)} \right\rfloor$$

We take

$$n = \frac{x}{2\log_2\left(\sqrt{a + x} + \sqrt{a}\right)}$$

We develop the equation to solve for $a$

$$\log_2\left(\sqrt{a + x} + \sqrt{a}\right) = \frac{x}{2n}$$

We apply $anti\log_2 N$ to both sides of the ecuation

$$\left(\sqrt{a + x} + \sqrt{a}\right) = 2^{x/2n}$$

Squaring both sides of the equation, we have

$$2a + x + 2\sqrt{a(a + x)} = 2^{x/n}$$

$$2\sqrt{a(a + x)} = 2^{x/n} - x - 2a$$

$$2\sqrt{a(a + x)} = \left(2^{x/n} - x\right) - 2a$$

Squaring both sides of the equation, we have

$$4a(a + x) = \left(2^{x/n} - x\right)^2 - 4a\left(2^{x/n} - x\right) + 4a^2$$

$$4a(a + x) + 4a\left(2^{x/n} - x\right) - 4a^2 = \left(2^{x/n} - x\right)^2$$

$$4a\left(a+x+2^{x/n}-x-a\right)=\left(2^{x/n}-x\right)^2$$

$$4a*2^{x/n}=\left(2^{x/n}-x\right)^2$$

$$a=\frac{\left(2^{x/n}-x\right)^2}{4*2^{x/n}}$$

Example. Find the value of $a$ when $n=8,\quad x=100$

We verify

$$\frac{100}{\log_2 100}>8>\frac{100}{2\log_2 100}$$

holds, then

$$a=\frac{\left(2^{100/8}-100\right)^2}{4*2^{100/8}}=1398$$

We have that between $1398$ and $1498$ there are at least $8$ primes numbers.

We verify by the proposition of the conjecture II. The lower bound for $\pi\left(1398,1498\right)$ is

$$\left\lfloor\frac{100}{2\log_2\left(\sqrt{1498}+\sqrt{1398}\right)}\right\rfloor=8$$

It is noteworthy that

$$\left(a+x\right)=\frac{\left(2^{x/n}+x\right)^2}{4*2^{x/n}}$$

We therefore conclude that between

$$\frac{\left(2^{x/n}-x\right)^2}{4*2^{x/n}}\quad\text{and}\quad\frac{\left(2^{x/n}+x\right)^2}{4*2^{x/n}}$$

there are at least $n$ primes numbers, for

$$\frac{x}{\log_2 x}\geq n\geq\frac{x}{2\log_2 x}$$

## 5.4. Algorithms for the calculation of $b$ given the values of $x$ and $n$.

**5.4.1.** If

$$\frac{x}{\log_2 x}\geq n$$

then

$$b=\frac{\left(x^{2}+2n.\log_{2}x\right)^{2}}{8n.x.\log_{2}x}$$

Obtaining the result.

By the conjecture III. we have that if $2\le b-a$, $\sqrt{b}-\sqrt{a}\le 1$, then the lower bund of $\pi(a,b)$ is

$$\left\lfloor\frac{(b-a)\left(\sqrt{b}-\sqrt{a}\right)^{2}}{2\log_{2}(b-a)}\right\rfloor$$

Doing $a=(b-x)$ we have that the lower bound for $\pi(a,b)$ is

$$\left\lfloor\frac{x\left(\sqrt{b}-\sqrt{b-x}\right)^{2}}{2\log_{2}x}\right\rfloor$$

We take

$$n=\frac{x\left(\sqrt{b}-\sqrt{b-x}\right)^{2}}{2\log_{2}x}$$

We develop the equation to solve for $b$ similar as we did it in **5.3.1**. until we get

$$b=\frac{\left(x^{2}+2n.\log_{2}x\right)^{2}}{8n.x.\log_{2}x}$$

Example. Find the value of $b$ when $x=200$, $n=3$

We verify

$$\frac{200}{2\log_{2}200}>n$$

holds, then

$$b=\frac{\left(200^{2}+6\log_{2}200\right)^{2}}{24*200*\log_{2}200}=43708$$

We have that between $43508$ and $43708$ there are at least $3$ prime numbers, according with the solution of the example **5.3.1.**

**5.4.2.** If

$$\frac{x}{\log_2 x} \geq n \geq \frac{x}{2\log_2 x}$$

then

$$b = \frac{\left(2^{x/n} + x\right)^2}{4*2^{x/n}}$$

Obtaining the result.

By the proposition of the conjecture II. we have that $4 \leq b-a$, $\sqrt{b}-\sqrt{a} \geq 1$, then the lower bound for $\pi\left(a,b\right)$ is

$$\left\lfloor \frac{b-a}{2\log_2\left(\sqrt{b}+\sqrt{a}\right)} \right\rfloor$$

Doing $a = \left(b-x\right)$ we have that the lower bound for $\pi\left(a,b\right)$ is:

$$\left\lfloor \frac{x}{2\log_2\left(\sqrt{b}+\sqrt{b-x}\right)} \right\rfloor$$

We take

$$n = \frac{x}{2\log_2\left(\sqrt{b}+\sqrt{b-x}\right)}$$

We develop the equation to solve for $b$ similar as we did it in **5.3.2**. until we get

$$b = \frac{\left(2^{x/n} + x\right)^2}{4*2^{x/n}}$$

Example. Find the value of $b$ when $n=8$, $x=100$

We verify

$$\frac{100}{\log_2 100} > 8 > \frac{100}{2\log_2 100}$$

holds, then

$$b = \frac{\left(2^{100/8} - 100\right)^2}{4*2^{100/8}} = 1498$$

We have that between $1398$ and $1498$ there are at least $8$ prime numbers, according with the solution of the example **5.3.2.**

We end this article with a famous quote about prime numbers:

> *"There are two facts about the distribution of prime numbers of which I hope to convince you so overwhelmingly that they will be permanently engraved in your hearts. The first is that, despite their simple definitions and role as the building blocks of the natural numbers, the prime numbers belong to the most arbitrary and ornery objects studied by mathematicians: they grow like weeds among the natural numbers, seeming to obey no other law than that of chance, and nobody can predict where the next one will sprout. The second fact is even more astonishing, for it states just the opposite: that the prime numbers exhibit stunning regularity, that there are laws governing their behavior, and that they obey these laws with almost military precision".* (DON ZAGIER, *1975)*

## 6. Conclusions.

If the conjectures I, II and III are true; then all problems related to the distance between prime numbers is solved, we have seen some problems known as The Andrica's conjecture are treated here by some of these conjectures, being reduced to cases special. The author believes that there is order in the sequence of prime numbers, their distribution is not chaotic, and is not governed by the laws of chance, and that given a prime number if we can predict where it will appear the next, it all depends on the validity of the conjectures discussed in this article.

The proof of these conjectures perhaps require advanced mathematical techniques, but his statements are simple and understanding is available to anyone with a basic knowledge of mathematics and the reader will seem that the question of the distribution of prime numbers is a pleasant and fun game. It is amazing that such simple results have not been previously discovered, when the ideas presented here could have been discovered a hundred years ago, and when hundreds of mathematicians (who squander talent and originality) have focused their efforts on this issue.

Perhaps the secret of the distribution of prime numbers is in this article. If the conjectures I, II and III are true; then an ancient and elusive problem will be solved forever.

**About the Author.**

*Name:* Jonas Castillo Toloza
*email*: jonas@fjcaldas.edu.co